\documentclass[3p]{elsarticle}

\usepackage{graphicx}
\usepackage{amsmath}
\usepackage{esint}
\usepackage{dirtytalk}
\usepackage[normalem]{ulem}
\usepackage{natbib}
\usepackage[font=scriptsize]{subfig}
\usepackage{epstopdf}
\usepackage{multirow}
\usepackage[colorlinks = true,linkcolor = blue,anchorcolor =red,citecolor = blue,filecolor = red,urlcolor = red,pdfauthor=author]{hyperref}
\usepackage[dvipsnames]{xcolor}
\usepackage{mwe}
\usepackage{makecell}
\usepackage{hhline}
\usepackage{float}
\usepackage{bookmark}
\usepackage{tabu}
\usepackage{relsize} 
\usepackage{empheq} 
\usepackage{upgreek} 
\usepackage{enumitem} 
\usepackage[bottom]{footmisc} 
\usepackage{amsfonts} 
\usepackage{derivative} 
\usepackage{colortbl} 
\usepackage{booktabs} 

\journal{}

\begin{document}
\interfootnotelinepenalty=10000

\begin{frontmatter}
	\title{Density-Based Topology Optimization of High-Fidelity Fluid-Structure Interaction Problems with Large Deformations}
	\author{Mohamed Abdelhamid}
	\author{Aleksander Czekanski\corref{mycorrespondingauthor}}
	\address{Department of Mechanical Engineering, Lassonde School of Engineering, York University, Toronto, Ontario, Canada}
	\cortext[mycorrespondingauthor]{Corresponding author}
	\ead{alex.czekanski@lassonde.yorku.ca}
	\begin{abstract}
		The application of modern topology optimization techniques to single physics systems has seen great advances in the last three decades. However, the application of these tools to sophisticated multiphysics systems such as fluid-structure interactions is still lagging behind, mainly due to the multidisciplinary and complex nature of such systems. In this work, we implement topology optimization of high-fidelity, fully-coupled fluid-structure interaction problems with large deformations. We use the arbitrary Lagrangian-Eulerian approach to deform the fluid mesh as a pseudo-structural system such that structural deformations are completely reflected in the fluid flow mesh. The fluid-structure interaction problem is formulated using the three-field formulation and the sensitivity analysis is derived using the discrete adjoint approach. We show through numerical examples the effect of the projection and interpolation parameters on the convergence and topology of the optimized designs and demonstrate the effect of considering the structural deformations in the fluid mesh.
	\end{abstract}
\end{frontmatter}


\section{Introduction}

The interaction of a deformable or movable structure with a flowing or a stationary fluid is a widespread phenomenon in nature known as \textbf{f}luid-\textbf{s}tructure \textbf{i}nteractions (FSI). Examples of FSI in everyday life include the flow of blood through flexible arteries, the deflection of an elastic airplane wing under aerodynamic loads, and the vibration of bridges and tall buildings under wind loads \cite{Richter2017}. Great advances have been reached in the numerical analysis of these problems \cite{Hou2012, Dowell2001}, and most multiphysics simulation packages - commercial and open source - include some degree of capability for solving these problems. Nonetheless, the application of numerical design tools such as \textbf{t}opology \textbf{o}ptimization (TO) to these problems is still lagging behind. This is mainly due to the multidisciplinary nature of the combined \textbf{t}opology \textbf{o}ptimization of \textbf{f}luid-\textbf{s}tructure \textbf{i}nteraction (TOFSI) problem and its strong nonlinear behavior, resulting in numerous stability and convergence issues.

A closer look at the chronology of TOFSI problems in literature adds valuable context to this discussion. The seminal work of \citet{Bendsoe1988} and \citet{Bendsoe1989} in the late 1980s marked the first appearance of modern density-based TO techniques, and most of the early literature that followed was dedicated to structural mechanics applications. In density-based TO, each discretization unit (e.g. a finite element) is interpolated between material and void to avoid the need for binary programming. To approach a discrete 0/1 state, penalization is enforced on a property of interest (e.g. Young's modulus in structural mechanics) to increase the cost of using intermediate densities \cite{Abdelhamid2021}. More than a decade later, the first application of TO to fluid flow problems appeared in the work of \citet{Borrvall2003} in 2003, where the authors introduced design parametrization through a parameter that controls permeability in Stokes fluid flow. Later in 2005, \citet{Gersborg-Hansen2005} extended the application of TO to Navier-Stokes equations using the analogy of a 2D channel with varying thickness for design parametrization. In 2005, \citet{Evgrafov2005} recognized the similarities between this design parametrization and Brinkman equations of fluid flow around the same time as \citet{Gersborg-Hansen2005}. Also of some relevance to TOFSI is the TO of structures under design-dependent loads where the fluid-structure interface is identified either explicitly \cite{Hammer2000, Fuchs2004} or implicitly \cite{Chen2001, Sigmund2007a}.

The first application of TO to FSI problems can be found in the pioneering work of \citet{Maute2002} in 2002 and \citet{Maute2004} in 2004 which presented the optimization of the internal structure of an aeroelastic wing. In these works, a three-field formulation was implemented to represent  the fluid, solid, and fluid mesh governing equations. Yet, the wet fluid-structure interface (i.e. the skin of the wing) was not included in the design domain, hence the term \textit{dry} TOFSI. Nonetheless, these works presented valuable contributions in demonstrating the effect of the coupled fluid-structure response on the optimized designs and in the sensitivity analysis of the coupled three-field system for TO purposes.

In contrast to dry TOFSI, \textit{wet} TOFSI includes the wet fluid-structure interface in the design domain. Optimizing this wet interface can be accomplished through one of two approaches: \textbf{(i)} moving the explicitly-defined fluid-structure interface during the optimization process, or \textbf{(ii)} interpolating between the fluid and solid within each discretization unit. The first approach can be accomplished through boundary tracking techniques such as the level-set method. The second approach necessitates overcoming a number of challenges for a successful implementation of the density-based approach to wet TOFSI. 

A few words on \textbf{the first approach} are in order. \citet{Jenkins2015, Jenkins2016} utilized the level-set method with the extended finite element method to perform dry TOFSI in \cite{Jenkins2015} and wet TOFSI in \cite{Jenkins2016}, both with large deformations. The level-set method, being a boundary tracking approach, inherently tracks the evolution of the explicit fluid-structure interface, hence the optimized design is usually of clear, crisp boundaries. However, the level-set approach lacks in the design freedom of internal features, which requires either modifying the level-set formulation to allow the creation of new holes, or seeding the initial design with a large number of internal holes. There is also attempts to reconcile the level-set method with density-based methods as in the work of \citet{Barrera2020}.

As for \textbf{the second approach}, it wasn't until the seminal work of \citet{Yoon2010} in 2010 that the first \textit{true} application of TO to FSI problems came to the light. \textit{True} TOFSI in the sense that both the internal features as well as the fluid-structure interface are included in the design domain, hence \textit{wet} TOFSI. Yoon's extension of \textit{dry} to \textit{wet} TOFSI using density-based methods can be attributed to three contributions:
\begin{enumerate}[label=\alph*.]
	\item The introduction of the \textit{unified} domain formulation where the solid and the fluid domains overlap such that solid and fluid coexist within the same discretization unit (i.e. finite element). This is in contrast to the traditional \textit{segregated} domain formulation where the solid and the fluid only ``meet" at the fluid-structure interface and the forces are transferred through surface stresses.
	\item The introduction of the force coupling filtering function\footnote{The filter function used to couple the fluid and the structure in TOFSI problems was first introduced in \citep[p.~599]{Yoon2010} and later termed ``pressure coupling filter function" in \citep[p.~972]{Lundgaard2018} as only pressure stresses were considered. In a previous manuscript \cite{Abdelhamid2022}, we extended the force coupling to total stress instead of the original pressure stress. Henceforth, we opted to generalize the term to be the ``force coupling filter function".} enabled the inclusion of low (or zero) intermediate density elements in the design domain without fearing severe mesh deformations as the amount of force transferred from the fluid to a solid element is directly related to its density value. This function is discussed in better detail in section \ref{ch7_sec:finite}.
	\item The transformation of the surface integral representation used in traditional FSI problems - where an explicit fluid-structural interface is defined - to a volumetric integral representation through the use of the divergence theorem. 
\end{enumerate}


Yoon followed his original work with a number of relevant articles such as stress-based TOFSI in \cite{Yoon2014}, design of compliant passive micro valves in \cite{Yoon2014a}, brittle and ductile failure of stress-based TOFSI in \cite{Yoon2017}, and transient TOFSI in \cite{Yoon2023}. Nonetheless, these works are limited to infinitesimal strain assumptions where the structural deformations are small and no mesh deformation of the fluid computational domain is implemented.

In this work, we extend the application of the density-based method to fully-coupled, wet TOFSI problems with large deformations. We implement the \textbf{A}rbitrary \textbf{L}agrangian-\textbf{E}ulerian approach to enable the deformation of the fluid mesh independently from the fluid motion. The rest of this work is organized as follows. In \textbf{section \ref{ch7_sec:gov_eqns}}, we discuss the governing equations of the TOFSI problem. In \textbf{section \ref{ch7_sec:finite}}, we detail the finite element formulations of the TOFSI problem as well as the sensitivity analysis utilized. In \textbf{section \ref{ch7_sec:dsgn_prblm}}, we describe the design problem to be optimized. In \textbf{section \ref{ch7_sec:numerc_exp}}, we discuss the setup of the numerical experiments and solve several cases to elaborate on the effect of the interpolation and projection parameters and the impact of considering vs neglecting mesh deformation. Finally, we conclude our work in \textbf{section \ref{ch7_sec:conclusion}}.

\section{Governing Equations of the Fluid-Structure Interaction Problem}
\label{ch7_sec:gov_eqns}

To enable the two-way coupling in an FSI problem with large deformations, we implement the \textit{three-field arbitrary Lagrangian-Eulerian formulation}, where the fluid mesh nodes may be fixed with the Eulerian description of the fluid, moved with the structural continuum in the Lagrangian description, or deformed arbitrarily in between \cite{Donea2017}. This arbitrary deformation of the fluid mesh nodes is typically chosen such that the mesh remains admissible (i.e. no overlapping or inverted elements).

In the unified domain formulation, the fluid computational domain $\mathit{\Omega}_f$ spans the entire computational domain $\mathit{\Omega}$. The solid computational domain $\mathit{\Omega}_s$ is fully contained within the fluid computational domain. The solid non-design domain $\mathit{\Omega}_{sn}$, if any, along with the design domain $\mathit{\Omega}_d$ constitute the solid domain. The fluid non-design domain $\mathit{\Omega}_{fn}$ is implicitly defined as the set resulting from subtracting the solid domain from the entire domain. The mesh deformation domain $\mathit{\Omega}_m$ is the same as the fluid non-design domain. These relations can be summarized as follows:
\begin{gather}
	\mathit{\Omega}_f := \mathit{\Omega}, \\
	\mathit{\Omega}_s = \mathit{\Omega}_d \cup \mathit{\Omega}_{sn}, \\
	\mathit{\Omega}_d \subseteq \mathit{\Omega}_s \subset \mathit{\Omega}_f, \\
	\mathit{\Omega}_{fn} := \mathit{\Omega}_f \backslash \mathit{\Omega}_s, \\
	\mathit{\Omega}_m = \mathit{\Omega}_{fn}.
\end{gather}

\begin{figure}[b!]
	\centering
	\includegraphics[width=0.6\textwidth]{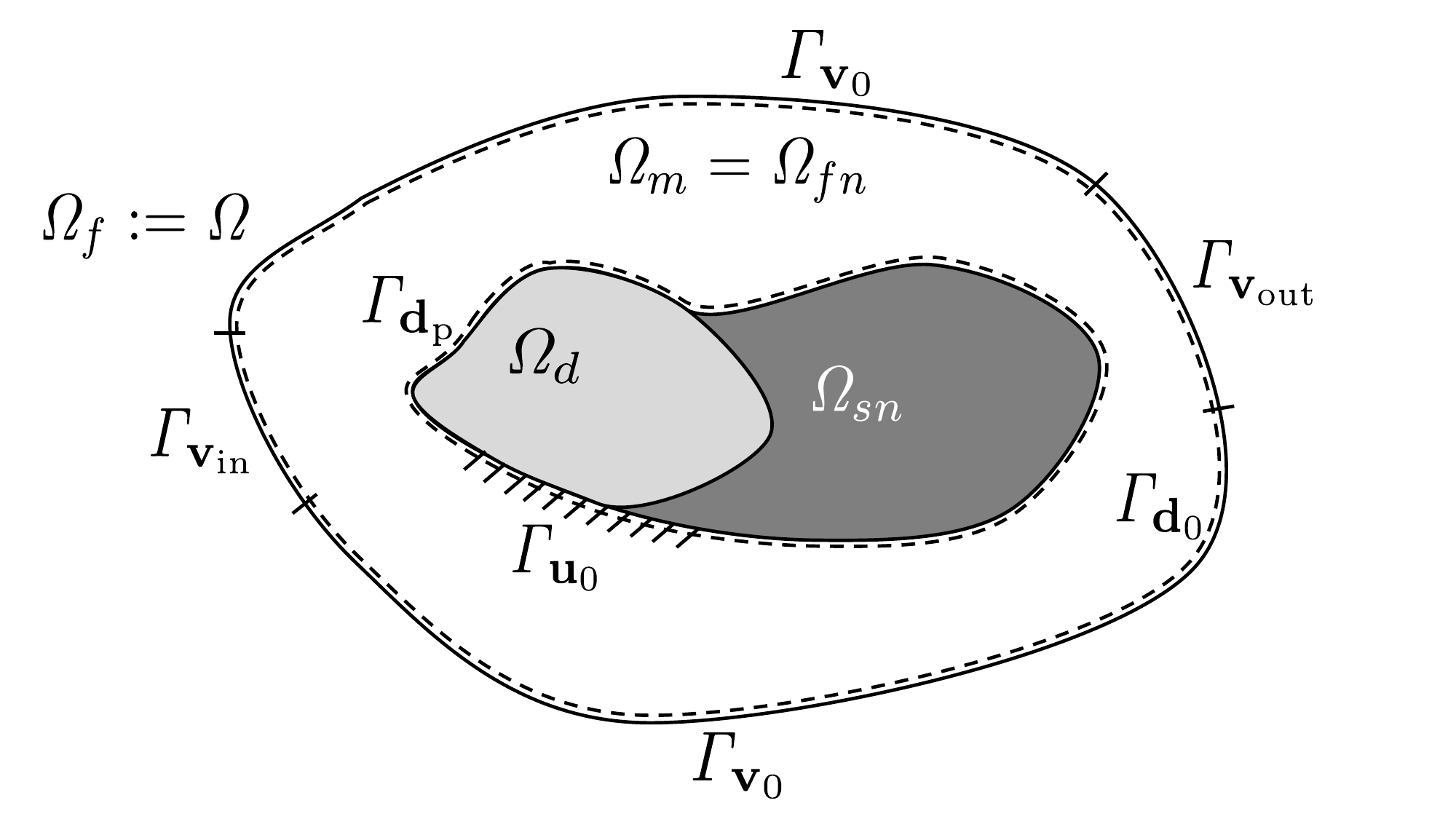}
	\caption{A schematic diagram of the TOFSI problem. Dotted lines indicate boundary identification for the mesh deformation domain.}
\end{figure}

In $\mathit{\Omega}_d$, a design variable $\rho$ per each finite element determines whether it's pure fluid, pure solid, or in between. More specifically, $\rho$ represents the solid density of the design element; 1 for pure solid, 0 for pure fluid, and porous solid in between. In density-based TO, the interpolation between fluid and solid is enabled through:
\begin{enumerate}[label=\alph*.]
	\item appending the fluid momentum equation with a Brinkman penalization term $-\alpha (\rho) \mathbf{v}$ (a scalar multiplied by the velocity vector) as a volume force where the inverse permeability\footnote{Yoon \cite[p.~502]{Yoon2010} called this term ``inverse permeability'' while Lundgaard et al. \cite[p.~971]{Lundgaard2018} called this term ``permeability''. We opted to use ``inverse permeability'' as this is in agreement with the original Brinkman equation \cite[p.~28]{Brinkman1947}.} $\alpha(\rho)$ is dependent on the design variable $\rho$, and 
	\item interpolating the structural stiffness $E(\rho)$, i.e. Young's Modulus, as a function of the design variable $\rho$.
\end{enumerate}

For a finite element that is 100\% solid (i.e. $\rho = 1$), the elastic modulus is effectively the bulk elastic modulus and the Brinkman penalization term enforces zero velocity within. On the other hand, for a finite element that is 100\% fluid (i.e. $\rho = 0$), the elastic modulus is effectively zero and the Brinkman penalization term is also zero such that the fluid velocity is unaffected. For a finite element with $0 < \rho < 1$, an interpolation between the two extreme states is applied. Therefore solid and fluid are allowed to co-exist within the same finite element, which is necessary to transform the discrete optimization problem into a continuous one.

In the following, we discuss the governing equations of the three-field formulation; the fluid flow, the structure, and the fluid mesh. For \textbf{the fluid flow}, the \textbf{N}avier-\textbf{S}tokes (NS) equations are used in their viscous, incompressible, steady-state form. The strong form of the partial differential equations is as follows \cite[p.~10]{Reddy2010}:
\begin{empheq}[right = \quad \empheqrbrace \, \text{in } \mathit{\Omega}_f]{gather}
	\boldsymbol{\nabla} \cdot \mathbf{v} = 0, \label{eq:contin} \\
	\rho_f \, (\mathbf{v} \cdot \boldsymbol{\nabla}) \mathbf{v} = \boldsymbol{\nabla} \cdot \boldsymbol{\upsigma}^f + \mathbf{f}^f - \alpha(\rho) \mathbf{v} , \label{ch7_eq:ns} \\ 
	\boldsymbol{\upsigma}^f = - p \mathbf{I} + \mu \left[ \boldsymbol{\nabla} \mathbf{v} + (\boldsymbol{\nabla} \mathbf{v})^T \right], \\
	\begin{split}
		\alpha(\rho) = \alpha_\text{max} + \left( 1 - \rho \right) \left( \alpha_\text{min} - \alpha_\text{max} \right) \frac{1 + p_\alpha}{1 - \rho + p_\alpha}. \label{ch7_eq:invrs_perm}
	\end{split}
\end{empheq}

\noindent where $\mathbf{v}$ is the fluid velocity, $\rho_f$ is the fluid density, $\boldsymbol{\upsigma}^f$ is the Cauchy fluid stress tensor, $\mathbf{f}^f$ is the external fluid force, $p$ is the hydrostatic pressure, and $\mu$ is the fluid dynamic viscosity. In Eq. \ref{ch7_eq:ns}, the last term on the right hand side is the Brinkman penalization term as discussed earlier. This form of the fluid momentum equation is henceforth designated the ``Brinkman-penalized NS" as opposed to the ``original NS" before adding the Brinkman penalization term. In Eq. \ref{ch7_eq:invrs_perm}, the inverse permeability $ \alpha (\rho) $ is dependent on the Brinkman penalization upper and lower limits, $\alpha_\text{max}$ and $\alpha_\text{min}$ respectively, and the Brinkman penalization interpolation parameter $p_\alpha$.

As for \textbf{the structure}, the Navier-Cauchy equations are used assuming linear elasticity with infinitesimal deformations under steady-state conditions. The strong form of the partial differential equations in the stress-divergence form is as follows \cite[p.~168]{Lai2010}:
\begin{empheq}[right = \quad \empheqrbrace \, \text{in } \mathit{\Omega}_s]{gather}
	\boldsymbol{\nabla} \cdot \boldsymbol{\upsigma}^s + \mathbf{f}^s = 0, \\
	\boldsymbol{\upsigma}^s = \mathbf{C}^s \, \boldsymbol{\upepsilon}^s, \\
	\boldsymbol{\upepsilon}^s = \frac{1}{2} \left[ \boldsymbol{\nabla} \mathbf{u} + (\boldsymbol{\nabla} \mathbf{u})^T \right], \\
	E(\rho) = E_\text{min} + (E_\text{max} - E_\text{min}) \rho^{\: \! p_E}. \label{ch7_eq:elastic_mod}
\end{empheq}

\noindent where $\boldsymbol{\upsigma}^s$ is the Cauchy solid stress tensor, $\mathbf{f}^s$ is the external solid force, $\mathbf{C}^s$ is the solid elasticity tensor under plane strain conditions\footnote{For FSI problems in 2D, the plane strain assumption makes more sense than the plane stress one since in the former the effects of the relatively large out-of-plane thickness may be ignored. In plane stress applications such as a microfluidic device, the effect of the relatively small out-of-plane thickness must be considered in the fluid flow equations (cf. \citep[p.~2]{abdelhamid2023calculation} and references therein).}, $\boldsymbol{\upepsilon}^s$ is the infinitesimal strain tensor, and $\mathbf{u}$ is the solid displacement. The elastic modulus $E(\rho)$ is interpolated using the modified SIMP approach where $E_\text{min}$ and $E_\text{max}$ are the lower and upper limits on $E$ representing the void and bulk elastic moduli, respectively, and $p_E$ is the elastic modulus penalization parameter.

As for the \textbf{fluid mesh}, it is assumed a pseudo-structural system with its own system of governing equations independent of the actual structure. Similarly to the structure, we employ the Navier-Cauchy equations assuming linear elasticity, infinitesimal deformations, and steady-state conditions. The strong form of the partial differential equations in the stress-divergence form is as follows \cite[p.~168]{Lai2010}:
\begin{empheq}[right = \quad \empheqrbrace \, \text{in } \mathit{\Omega}_m]{gather}
	\boldsymbol{\nabla} \cdot \boldsymbol{\upsigma}^m + \mathbf{f}^m = 0, \\
	\boldsymbol{\upsigma}^m = \mathbf{C}^m \, \boldsymbol{\upepsilon}^m, \\
	\boldsymbol{\upepsilon}^m = \frac{1}{2} \left[ \boldsymbol{\nabla} \mathbf{d} + (\boldsymbol{\nabla} \mathbf{d})^T \right].
\end{empheq}

\noindent where $\boldsymbol{\upsigma}^m$, $\mathbf{f}^m$, $\mathbf{C}^m$, and $\boldsymbol{\upepsilon}^m$ are similarly defined for the mesh pseudo-structure and  $\mathbf{d}$ is the fluid mesh displacement. The elastic modulus of the mesh pseudo-structure can be taken as any value where the only concern is the appropriate numerical scaling of the system. In this work, it's taken as unity.

Next, we move on to define the boundary conditions. The essential boundary conditions are defined as follows:
\begin{alignat}{3}
	\text{Fluid No-slip:} & \qquad \mathbf{v} = \mathbf{0} & \quad \text{on } & \mathit{\Gamma}_{\mathbf{v}_0}, \label{ch7_eq:bc_noslip} \\
	\text{Fluid Inlet:} & \qquad \mathbf{v} = \mathbf{v}_\text{in} & \quad \text{on } & \mathit{\Gamma}_{\mathbf{v}_\text{in}}, \label{eq:bc_v_in} \\
	\text{Fluid Outlet:} & \qquad p = 0 & \quad \text{on } & \mathit{\Gamma}_{\mathbf{v}_\text{out}}, \label{eq:bc_p_out} \\
	\text{Structural Zero Displacement:} & \qquad \mathbf{u} = \mathbf{0} & \quad \text{on } & \mathit{\Gamma}_{\mathbf{u}_0}, \\
	\text{Mesh Prescribed Displacement:} & \qquad \mathbf{d} = \mathbf{u} & \quad \text{on } & \mathit{\Gamma}_{\mathbf{d}_p}, \\
	\text{Mesh Zero Displacement:} & \qquad \mathbf{d} = \mathbf{0} & \quad \text{on } & \mathit{\Gamma}_{\mathbf{d}_0}.
\end{alignat}

\noindent where the boundaries are defined as:
\begin{gather}
	\mathit{\Gamma}_{\mathbf{v}_0} := \mathit{\partial \Omega}_f \backslash (\mathit{\Gamma}_{\mathbf{v}_\text{in}} \cup \mathit{\Gamma}_{\mathbf{v}_\text{out}}), \\
	\mathit{\Gamma}_{\mathbf{d}_0} := \mathit{\partial \Omega}_f = \mathit{\Gamma}_{\mathbf{v}_0} \cup \mathit{\Gamma}_{\mathbf{v}_\text{in}} \cup \mathit{\Gamma}_{\mathbf{v}_\text{out}}, \\
	\mathit{\Gamma}_{\mathbf{d}_p} := \mathit{\Omega}_{fn} \cap \mathit{\Omega}_s.
\end{gather}

Note that the fluid no-slip boundary condition in Eq. \ref{ch7_eq:bc_noslip} is only defined on $\mathit{\Gamma}_{\mathbf{v}_0}$, which is the external domain boundary aside from the inlet and outlet. The volume force term appended to the fluid momentum in Eq. \ref{ch7_eq:ns} automatically enforces a no-slip condition wherever needed within the solid domain and/or its fluid-structure interface. The natural boundary condition in an FSI problem with an explicitly-defined fluid-structure interface is typically defined as follows:
\begin{equation}
	\boldsymbol{\upsigma}^f \cdot \mathbf{n}^f = \boldsymbol{\upsigma}^s \cdot \mathbf{n}^s \qquad \text{on } \mathit{\Gamma}_\text{FSI}
	\label{eq:trac_eq}
\end{equation}

\noindent where $\mathbf{n}^f$ and $\mathbf{n}^s$ are the normals to the fluid and solid surfaces, respectively. While the fluid-structure interface is explicitly defined in pure FSI problems (i.e. pure is in no porous media), this is not the case in density-based TOFSI problems where all solids are porous to some degree. Hence, to apply the traction coupling condition, $\mathit{\Gamma}_\text{FSI}$ must be taken as all edges of all solid finite elements. Alternatively, the divergence theorem may be used to transform the surface integral into a volume integral as follows \citep[p.~598]{Yoon2010}:
\begin{equation}
	\mathbf{f}^s = \oiint\limits_{\mathit{\Gamma}_{FSI}} \boldsymbol{\upsigma}^f \cdot \mathbf{n}^f \odif{\mathit{\Gamma}} = \iiint\limits_{\mathit{\Omega}_s} \nabla \cdot \boldsymbol{\upsigma}^f \odif{\mathit{\Omega}}.
\end{equation}

In the next section, we transform the strong form of the governing equations into the discretized finite element form.

\section{Finite Element Formulations}
\label{ch7_sec:finite}

\subsection{Three-Field Formulation}

In pure FSI problems, conformal or matching discretization is typically used to ensure simple, yet accurate force coupling between the fluid and the structure. By conformal discretization we mean the fluid mesh nodes and the solid mesh nodes coincide at all times. While this conformal discretization is only needed at the explicitly-defined fluid-structure interface in pure FSI problems, this is not the case in density-based TOFSI problems. Instead, the fluid and solid nodes must match in the entire solid computational domain.

To represent the FSI problem, we utilize the three-field formulation originally proposed by \citet{Farhat1995} and later used in an optimization context in \cite{Maute2003}. In the discretized form, this formulation leads to the following set of governing equations:
\begin{equation}
	\mathbb{R}(\mathbf{\hat{r}}) = \left\{
	\begin{array}{l}
		\mathbb{S}(\boldsymbol{\uprho}, \mathbf{\hat{u}}, \mathbf{\hat{d}}, \mathbf{\hat{w}}) \\
		\mathbb{D}(\boldsymbol{\uprho}, \mathbf{\hat{u}}, \mathbf{\hat{d}}) \\
		\mathbb{F}(\boldsymbol{\uprho}, \mathbf{\hat{d}}, \mathbf{\hat{w}})
	\end{array} \right\} = \mathbf{0}
\end{equation}
\noindent where $\mathbb{S}$, $\mathbb{D}$, and $\mathbb{F}$ are the structural, fluid mesh, and fluid flow governing equations, respectively. $\mathbf{\hat{r}}$ is a vector of all state variables in the FSI problem (structural displacements $\mathbf{\hat{u}}$, fluid mesh deformations $\mathbf{\hat{d}}$, and fluid flow velocities and pressures $\mathbf{\hat{w}}$) in the discretized form. $\boldsymbol{\uprho}$ is a vector of all design variables.

For \textbf{the fluid flow}, the discretized finite element form of the NS equations is obtained through multiplying the strong form of the PDEs with suitable weight functions then integrating over the appropriate computational domain. Integration by parts is typically utilized to reduce the continuity requirements on the shape/interpolation functions used to approximate the solution through moving the differentiation to the weight functions \cite{Reddy2010}. We implement the \textit{standard Galerkin method of weighted residuals} where the shape and weight functions are the same. The resulting finite element model is of the mixed (i.e. velocity-pressure) type where the velocities and pressures are both unknown and solved for simultaneously. We implement a \textit{quadrilateral meshing} that is regular and structured wherever possible. To satisfy the \textit{Ladyzhenskaya-Babuska-Brezzi} condition \citep[p.~176]{Reddy2010}, we use \textit{Q2Q1} Lagrangian finite elements, that is 9 velocity nodes and 4 pressure nodes. At low to moderate Reynolds numbers, this finite element eleminates the need for implementing stabilization techniques that are typically used to dampen the discontinuities. The fluidic boundary conditions are applied strongly, i.e. node-wise, hence there are no externally applied fluid forces. The resulting nonlinear system is solved using the undamped Newton-Raphson method \citep[p.~190]{Reddy2010}.

In a previous work \cite{Abdelhamid2022}, we extended the force coupling in TOFSI from hydrostatic coupling to total stress coupling, where the main challenge was overcoming the discontinuities in the velocity derivatives at the elemental boundaries. In this work, we implement the total stress coupling developed earlier as it results a higher fidelity in the numerical representation of the FSI problem and avoids the issues associated with the pressure spikes in the hydrostatic force coupling as is discussed further in the numerical examples. After calculating the fluidic forces per each finite element and before the global assembly of the discretized force vector, the forces must be multiplied by the force coupling filtering function defined as \citep[p.~607]{Yoon2010}:
\begin{equation}
	\mathit{\Upsilon}(\rho_i) = \mathit{\Upsilon}_\text{min} + (\mathit{\Upsilon}_\text{max} - \mathit{\Upsilon}_\text{min}) \ \rho_i^{\, p_\mathit{\Upsilon}}
\end{equation}

\noindent where $\mathit{\Upsilon}_\text{max}$ and $\mathit{\Upsilon}_\text{min}$ are assumed 1 and 0 respectively. $p_\mathit{\Upsilon}$ is the force coupling filter function interpolation parameter. The main intention of using the force coupling filtering function is to prevent severe deformation of low density elements which might lead to potential singularities in solving the structural finite element system.

As for \textbf{the structure}, the Navier-Cauchy equation is transformed into the discretized finite element form through the \textit{principle of virtual displacement} following the procedure in \citep[p.~153]{Bathe2014finite}. The only deviation from that procedure is that the elastic modulus in the constitutive matrix of each finite element is taken as dependent on its density according to Eq. \ref{ch7_eq:elastic_mod}. To enable conformal meshing and to simplify the force coupling procedure, \textit{Quadratic Lagrangian} finite elements (i.e. 9 displacement nodes) are used for the structure. Similarly to the structure, \textbf{the fluid mesh} governing equations are discretized using the principle of virtual displacements and Quadratic Lagrangian finite elements are used \cite{Abdelh}. There are no external forces applied to the fluid mesh system as the boundary conditions are applied in a strong, node-wise manner. For the entire computational domain, all finite elements are assumed isoparametric, hence the geometry is interpolated through the same shape functions used to interpolate the state variables. The entire computational domain is discretized into 5050 finite elements in a regular, structured mesh.

\begin{figure}[!]
	\centering
	\includegraphics[height=0.95\textheight]{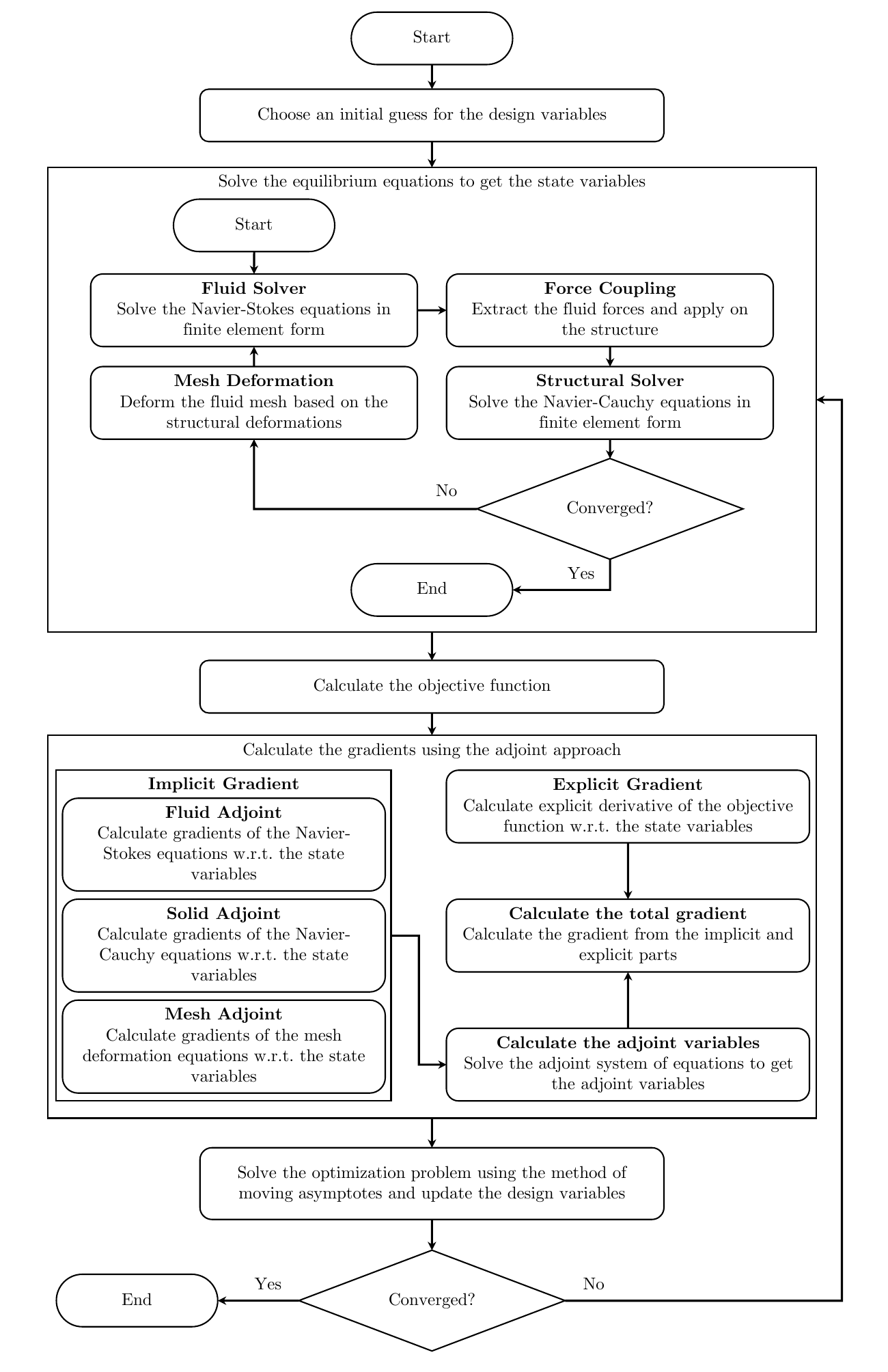}
	\caption{Flow chart of the entire TOFSI problem.}
	\label{ch7_fig:flowchart}
\end{figure}

One of the major difficulties resulting from considering large deformations in the TOFSI problem is the increased nonlinearity of the analysis problem requiring significant computational resources.  A flow chart of the entire TOFSI problem is presented in Fig. \ref{ch7_fig:flowchart} showing the numerical procedure followed for the analysis and optimization steps. For the analysis step, there are mainly two approaches for solving strongly coupled FSI systems; monolithic/simultaneous or segregated/staggered/partitioned \cite{Rugonyi2001}. Both terms refer to how the mathematical equations are handled during the solution of the governing equations and has nothing to do with the strength of the physical coupling between the fluid and the structure. An interesting discussion on the advantages and disadvantages of either approach is available in \citep[p.~916]{Maute2003}. Although some commercial software packages use the term ``fully coupled" to describe solving the entire system of nonlinear equations simultaneously, it can be misleading as to mean the strength of the physical coupling between the differing physics (e.g. fluid and structure) of the problem \citep[p.~1553]{comsol6}. We prefer using either ``simultaneous" or ``monolithic" as these terms refer directly to the way the mathematical equations are handled without ambiguity. On the other side, there is the ``staggered" or ``partitioned" solution procedure where the different physics are solved in steps with coupling conditions enforced from each step to the next. In this work, we implement the partitioned solution approach where each system of governing equations is solved on its own before moving to the next system. Convergence is achieved when the residuals of all the unknown state variables is lower than a certain tolerance.

In the next subsection, we discuss the derivation of the sensitivity analysis of the TOFSI problem using the adjoint method.

\subsection{Sensitivity Analysis of the Topology Optimization Problem}
\label{ch7_ssec:snst_anls}

Typically, in TO problems, the number of design variables greatly exceeds the number of constraints. Hence, it is better from a computational perspective to use the \textit{adjoint} method in calculating the sensitivities as opposed to the \textit{direct} method (cf. \cite{Maute2003} for a discussion on the adjoint method in an FSI context). Given an objective function (or a constraint) $f$, the sensitivity of $f$ w.r.t. the design variable $\rho_i$ is calculated as follows:
\begingroup
\allowdisplaybreaks
\begin{align}
	& \odv{f}{\rho_i} = \underbrace{\pdv{f}{\rho_i}}_{\text{explicit}} +	\underbrace{\pdv{f}{\mathbf{\hat{r}}} ^T \pdv{\mathbf{\hat{r}}}{\rho_i}}_{\text{implicit}}, \\
	& \odv{\mathbb{R}}{\rho_i} = \mathbf{0}, \quad \rightarrow \quad \pdv{\mathbb{R}}{\rho_i} + \pdv{\mathbb{R}}{\mathbf{\hat{r}}} \pdv{\mathbf{\hat{r}}}{\rho_i} = \mathbf{0}, \\
	& \pdv{\mathbf{\hat{r}}}{\rho_i} = - \pdv{\mathbb{R}}{\mathbf{\hat{r}}} ^{-1}
	\pdv{\mathbb{R}}{\rho_i}, \\
	& \odv{f}{\rho_i} = \pdv{f}{\rho_i} -	\pdv{f}{\mathbf{\hat{r}}} ^T \pdv{\mathbb{R}}{\mathbf{\hat{r}}} ^{-1} \pdv{\mathbb{R}}{\rho_i}, \\
	& \odv{f}{\rho_i} = \pdv{f}{\rho_i} -
	\left( \pdv{\mathbb{R}}{\mathbf{\hat{r}}} ^{-T}	\pdv{f}{\mathbf{\hat{r}}} \right)^T		\pdv{\mathbb{R}}{\rho_i}, \\
	& \boldsymbol{\uplambda} = \pdv{\mathbb{R}}{\mathbf{\hat{r}}} ^{-T}
	\pdv{f}{\mathbf{\hat{r}}}, \quad \rightarrow \quad \pdv{\mathbb{R}}{\mathbf{\hat{r}}} ^T \boldsymbol{\uplambda} = \pdv{f}{\mathbf{\hat{r}}}, \\
	& \odv{f}{\rho_i} = \pdv{f}{\rho_i} - \boldsymbol{\uplambda}^T
	\pdv{\mathbb{R}}{\rho_i}.
\end{align}
\endgroup

\noindent where $\boldsymbol{\uplambda}$ is the adjoint vector which is distinctive for each objective function or constraint.

\section{Description of the Design Problem}
\label{ch7_sec:dsgn_prblm}

The \textit{original} version of the column in a channel test problem was first discussed in a TOFSI context in \cite[p.~610]{Yoon2010} and has been used later as a benchmark problem in a number of works on TOFSI. It was later modified by \citet{Lundgaard2018} to increase the relative size of the design domain w.r.t. the entire computational domain, rescale it from the micro to the macro scale, and generally strengthen the fluid-structure dependency. In this work, we use the modified version as the test problem, henceforth the designation ``modified" is dropped.

\begin{figure}[b!]
	\centering
	\includegraphics[width=0.6\textwidth]{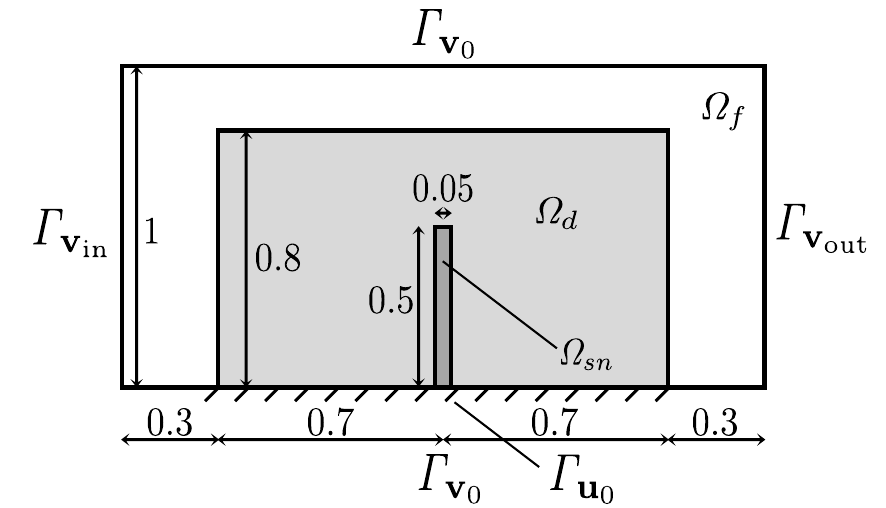}
	\caption{The \textit{modified} column in a channel design problem as described in \cite{Lundgaard2018}.}
	\label{ch7_fig:mdfd_clmn_in_a_chnl}
\end{figure}

The problem consists of a 0.8 $\times$ 1.4 m rectangular design space (light gray) placed inside a 1 $\times$ 2 m rectangular channel (cf. Fig. \ref{ch7_fig:mdfd_clmn_in_a_chnl}). A 0.05 $\times$ 0.5 m non-design elastic column (dark gray) is placed within the design space to force the optimizer to reach a more sophisticated solution than a simple bump at the bottom of the channel. The top and bottom surfaces of the channel $\mathit{\Gamma}_{\mathbf{v}_0}$ have a no-slip condition applied. A fully-developed, parabolic laminar flow profile is applied at the inlet $\mathit{\Gamma}_{\mathbf{v}_\text{in}}$ on the left and a zero pressure condition is applied at the outlet $\mathit{\Gamma}_{\mathbf{v}_\text{out}}$ on the right. The bottom surface of the design and non-design spaces $\mathit{\Gamma}_{\mathbf{u}_0}$ is fixed to a ground structure.

Unless otherwise noted, the following default values are used:
\begin{gather}
	v_\mathrm{max} = 1 \ \mathrm{m}/\mathrm{s}, \\
	\rho_f = 1 \ \mathrm{kg}/\mathrm{m}^3, \\
	\mu = 1 \ \mathrm{Pa} \cdot \mathrm{s}, \\
	E_\mathrm{max} = 1e+4 \ \mathrm{Pa}, \\
	E_\mathrm{min} = 1e-6 \ \mathrm{Pa}, \\
	\nu = 0.3.
\end{gather}

\noindent where $v_\text{max}$ is the maximum velocity at the parabolic inlet profile $\mathit{\Gamma}_{\mathbf{v}_\text{in}}$ and $\nu$ is the structural Poisson's ratio. The minimum limit on the elastic modulus $E_\text{min}$ is taken as $1e-10$ of $E_\text{max}$ and not zero to avoid singularities since, in the modified SIMP approach, the design variables are allowed to be zero.

The initial guess is taken as follows; the usual volume fraction $v_f$ for elements in the design domain, 1 for elements in the solid non-design domain, and 0 for elements in the fluid non-design domain as detailed in Eq. \ref{ch7_eq:init_guss}:
\begin{empheq}[left = {\rho = \empheqlbrace}]{equation}
	\begin{aligned}
		& v_f & \mathrm{in} \quad & \Omega_d, \\
		& 1 & \mathrm{in} \quad & \Omega_{sn}, \\
		& 0 & \mathrm{in} \quad & \Omega_{fn}.
	\end{aligned}
	\label{ch7_eq:init_guss}
\end{empheq}

\section{Setup and Results of the Numerical Experiments}
\label{ch7_sec:numerc_exp}

In this section, we discuss the TO problem definition and the setup used for the numerical experiments. The design problem, described in section \ref{ch7_sec:dsgn_prblm}, is solved with the objective function of minimizing the compliance of the structure at a specific volume fraction as follows:
\begin{equation}
	\begin{aligned}
		\text{minimize:} & \quad f = \sum_i \mathbf{\hat{u}}_i ^T \ \mathbf{k}_i^s \ \mathbf{\hat{u}}_i, \\[4pt]
		\text{subject to:} & \quad \mathbf{0} \leq \boldsymbol{\uprho} \leq \boldsymbol{1} , \\[4pt]
		& \quad \sum_i V_i/V_0 \leq V_f .
	\end{aligned}
	\label{ch7_eq:to_prb_dfntn}
\end{equation}

\noindent where $\mathbf{\hat{u}}_i$ and $\mathbf{k}_i^s$ are the nodal displacement vector and stiffness matrix of finite element $i$, $V_i$ is the volume of finite element $i$, $V_0$ is the volume of the entire design domain, and the volume fraction $V_f$ is set to 0.1 in agreement with \citep[p.~613]{Yoon2010} and \citep[p.~976]{Lundgaard2018}. The summation in Eq. \ref{ch7_eq:to_prb_dfntn} spans all elements in the solid computational domain, both design and non-design. Recall that the volume fraction is further used as the initial guess for finite elements in the design domain as indicated earlier. The entire problem is coded and solved in MATLAB.

The problem is solved using a combination of robust formulation and non-unity $p_E$ and $p_\mathit{\Upsilon}$. For the robust formulation, we followed Eqs. 10-13 in \citep[p.~973]{Lundgaard2018} with projection thresholding values of $\eta_n = 0.50$, $\eta_d = 0.49$, and $\eta_e = 0.51$ for the nominal, dilated, and eroded designs, respectively. The strength of the projection in the robust formulation is determined by a parameter $\beta$ where more discrete designs are obtained through increasing $\beta$. The relative filtering radius is set to 1.5 which is the smallest possible value for preventing checkerboarding and enabling the use of the robust formulation.

Continuation is implemented for $\beta$ as well as $p_E$ and $p_\mathit{\Upsilon}$. The problem is started with $p_E = 1$, $p_\mathit{\Upsilon} = 1$, and $\beta = 4$. These parameters are updated such that at iterations \{21, 41, 61, 81\}, $\beta$ is set to \{8, 16, 32, 64\}, while $p_E$ and $p_\mathit{\Upsilon}$ are raised by 0.5 and 0.5/$\delta_{E|\mathit{\Upsilon}}$, respectively. $\delta_{E|\mathit{\Upsilon}}$ is a parameter that is used to modify the penalization of $p_\mathit{\Upsilon}$ against that of $p_E$. For $\delta_{E|\mathit{\Upsilon}} = 1$, $p_E$ and $p_\mathit{\Upsilon}$ experience the same degree of penalization. As $\delta_{E|\mathit{\Upsilon}}$ increases (or decreases), $p_\mathit{\Upsilon}$ gets less (or more) penalized than $p_E$. The total number of iterations is 100 such that each continuation step uses 20 iterations.

%

The \textbf{m}ethod of \textbf{m}oving \textbf{a}symptotes (MMA) is used to update the design variables \cite{Svanberg1987}. We used the default settings of the MMA code in its min/max form \citep[p.~3]{Svanberg2004} except for the following changes: \textbf{(i)} the move limit is set to 0.1, and \textbf{(ii)} a positive offset of 1 is added to the objective functions before calling the MMA function to increase the constraint violation. These modifications increase the ``aggressiveness" of the MMA code such that, during the optimization process, the maximum change in the design variables is typically above 0.05 and even higher at the start and after each continuation update. The volume fraction constraint is applied on the dilated design and its value is updated every two iterations to account for this increased aggressiveness, cf. Eq. 14 in \citep[p.~774]{Wang2011} for more details.

\subsection{Verification of the Sensitivity Analysis}

In this subsection, we verify the analytically-derived sensitivity analysis in subsection \ref{ch7_ssec:snst_anls} using the complex-step derivative approximation method \cite{Martins2003}. This method is implemented through adding a small imaginary step $ih$ to a single design variable $\rho_i$ at a time, and calculating the approximate derivative of the objective function (or constraint) $f$ w.r.t. that particular design variable $\rho_i$ as follows:
\begin{equation}
	\pdv{f}{\rho_i} \approx \frac{\text{Im} \left[ f(\rho_i + ih) \right] }{h} 
\end{equation} 

\noindent where $\text{Im}(f)$ indicates the imaginary part of $f$.

\begin{figure}[t!]
	\centering
	\includegraphics[width=0.5\textwidth]{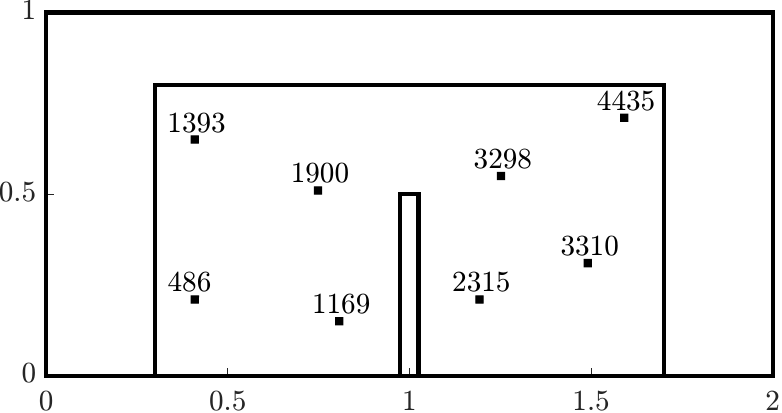}
	\caption{Location and number designations of the finite elements selected for sensitivity verification.}
	\label{ch7_fig:lctn_of_snst_elmnts}
\end{figure}

For the column in a channel test problem described in section \ref{ch7_sec:dsgn_prblm}, eight finite elements are sampled within the design domain. For these elements, the sensitivities are calculated using the analytically-derived equations in subsection \ref{ch7_ssec:snst_anls} and compared to those calculated through the complex-step derivative approximation method. Figure \ref{ch7_fig:lctn_of_snst_elmnts} shows the location of the finite elements of interest within the design domain. For the results in table \ref{tab:snst_vrfctn}, we utilized $p_\alpha = 18e-7$, $p_E = 1$, $p_\Upsilon = 1$, and $\beta = 4$ in addition to the initial guess described in Eq. \ref{ch7_eq:init_guss}. The errors reported in table \ref{tab:snst_vrfctn} are extracted for the nominal design and are calculated by subtracting the analytical sensitivity from the complex sensitivity then dividing over the analytical sensitivity. As can be noted, the errors are within a reasonable range attesting to the validity of the derived analytical sensitivity. It's worth noting that due to the strong nonlinearity of the TOFSI problem with mesh deformation, numerical errors tend to compound if the convergence tolerance is not tight enough. From our experience, a convergence tolerance of $1e-3$ may be appropriate for TOFSI problems without mesh deformation, but produced significant errors in the sensitivity verification of TOFSI problems with mesh deformation, hence we had to implement a much tighter tolerance of $1e-8$.

\begin{table}[t!]
	\centering
	\captionsetup{width=0.9\textwidth}
	\caption{Verification of analytically-derived sensitivities using the complex-step derivative approximation method.}
	\label{tab:snst_vrfctn}
	\def\arraystretch{1.5}
	\begin{tabular}{@{}crrrr@{}} \toprule
		Element & Imaginary Step & Obj. Fun. Imaginary Part & Analytical Sensitivity & Normalized Error \\ \midrule
		1393 & $+1e-10$ & $+8.289833322905e-13$ & $+8.289833324822e-03$ & $-2.3124e-08\%$ \\
		486 & $+1e-10$ & $-1.180951947628e-12$ & $-1.180951947740e-02$ & $-9.4322e-09\%$ \\
		1900 & $-1e-10$ & $-8.269116672568e-13$ & $+8.269116670589e-03$ & $+2.3934e-08\%$ \\
		1169 & $-1e-10$ & $+3.847814031866e-12$ & $-3.847814031741e-02$ & $+3.2446e-09\%$ \\
		3298 & $+1e-10$ & $+1.746191997720e-12$ & $+1.746191998103e-02$ & $-2.1917e-08\%$ \\
		2315 & $+1e-10$ & $-3.112693493737e-12$ & $-3.112693493588E-02$ & $+4.8049e-09\%$ \\
		4435 & $-1e-10$ & $+8.989216991614e-13$ & $+8.989216995844e-03$ & $-4.7067e-08\%$ \\
		3310 & $+1e-10$ & $-9.146870951814e-13$ & $-9.146870950669e-03$ & $+1.2521e-08\%$ \\
		\bottomrule
	\end{tabular}
\end{table}

\subsection{On the Selection of the Projection and Interpolation Parameters}

\begin{figure}[h!]
	\centering
	\captionsetup[subfigure]{justification=centering}
	\subfloat[$p_\alpha = 8.0e-7$, $f = 0.153763$, \\ \vspace{1pt} $DM = 0.1373 \ \%$]{\includegraphics[width=0.32\textwidth]{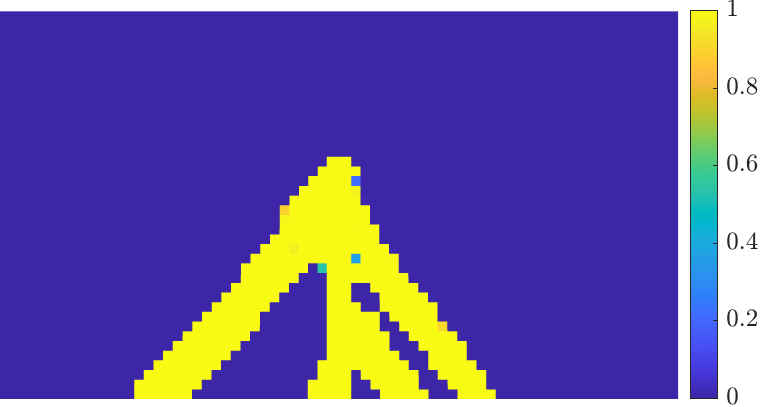}} \
	\subfloat[$p_\alpha = 8.5e-7$, $f = 0.151441$, \\ \vspace{1pt} $DM = 0.044006 \ \%$]{\includegraphics[width=0.32\textwidth]{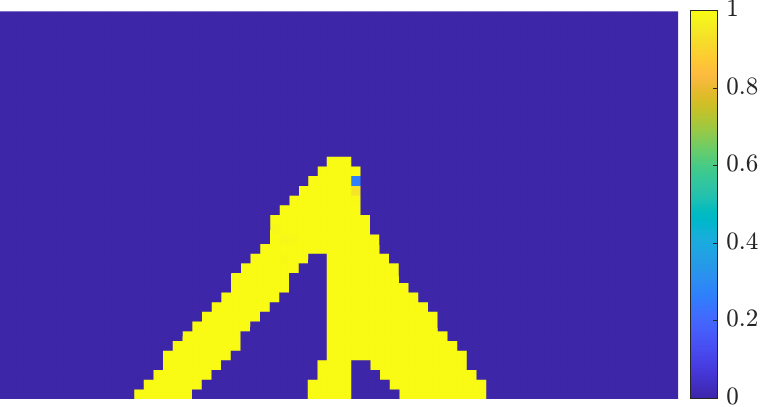}} \
	\subfloat[$p_\alpha = 9.0e-7$, $f = 0.151564$, \\ \vspace{1pt} $DM = 0.059782 \ \%$]{\includegraphics[width=0.32\textwidth]{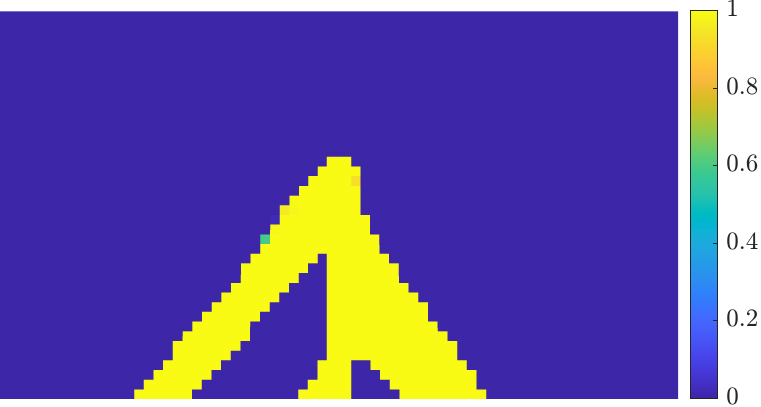}} \\
	\subfloat[$p_\alpha = 9.5e-7$, $f = 0.150005$, \\ \vspace{1pt} $DM = 0.047897 \ \%$]{\includegraphics[width=0.32\textwidth]{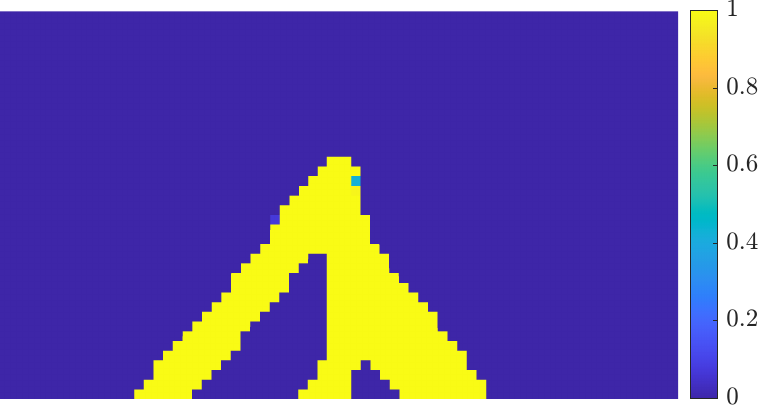}} \
	\subfloat[$p_\alpha = 10.0e-7$, $f = 0.152509$, \\ \vspace{1pt} $DM = 0.13620 \ \%$]{\includegraphics[width=0.32\textwidth]{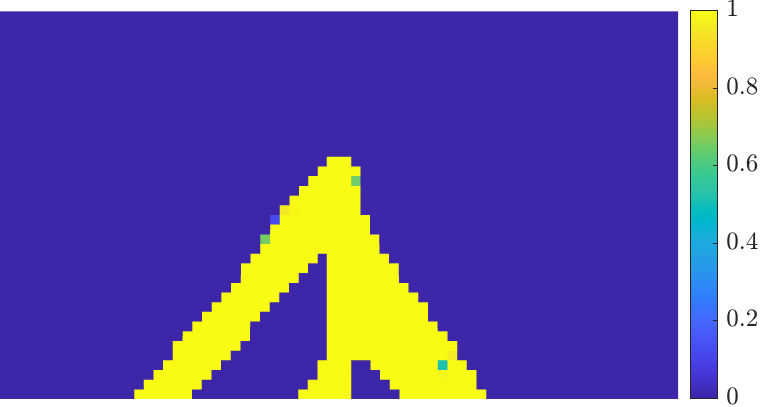}} \
	\subfloat[$p_\alpha = 11.5e-7$, $f = 0.152476$, \\ \vspace{1pt} $DM = 0.031296 \ \%$]{\includegraphics[width=0.32\textwidth]{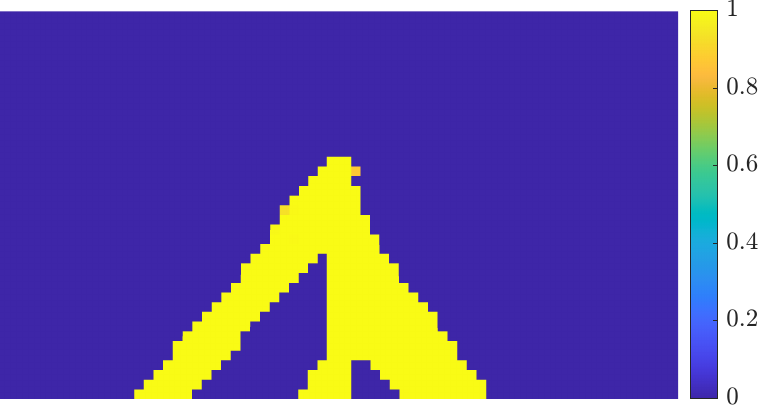}} \\
	\subfloat[$p_\alpha = 12.5e-7$, $f = 0.153079$, \\ \vspace{1pt} $DM = 0.075258 \ \%$]{\includegraphics[width=0.32\textwidth]{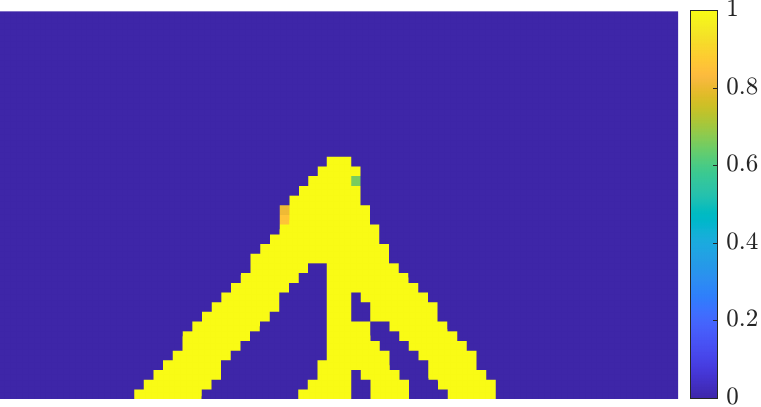}} \
	\subfloat[$p_\alpha = 14.0e-7$, $f = 0.159349$, \\ \vspace{1pt} $DM = 0.019760 \ \%$]{\includegraphics[width=0.32\textwidth]{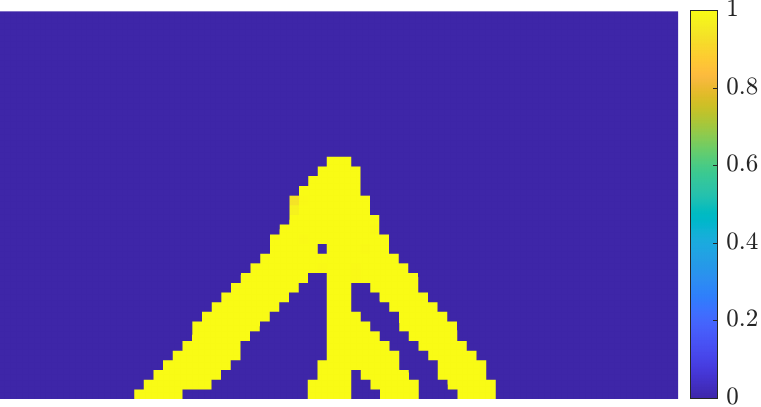}} \
	\subfloat[$p_\alpha = 15.5e-7$, $f = 0.151417$, \\ \vspace{1pt} $DM = 0.074926 \ \%$]{\includegraphics[width=0.32\textwidth]{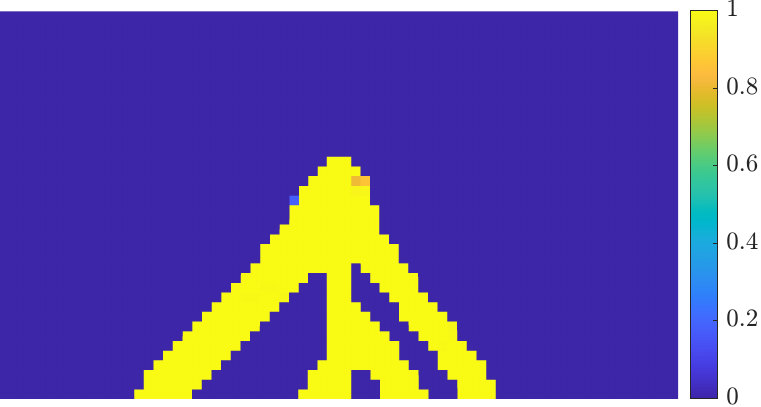}} \\
	\subfloat[$p_\alpha = 17.5e-7$, $f = 0.152704$, \\ \vspace{1pt} $DM = 0.080990 \ \%$]{\includegraphics[width=0.32\textwidth]{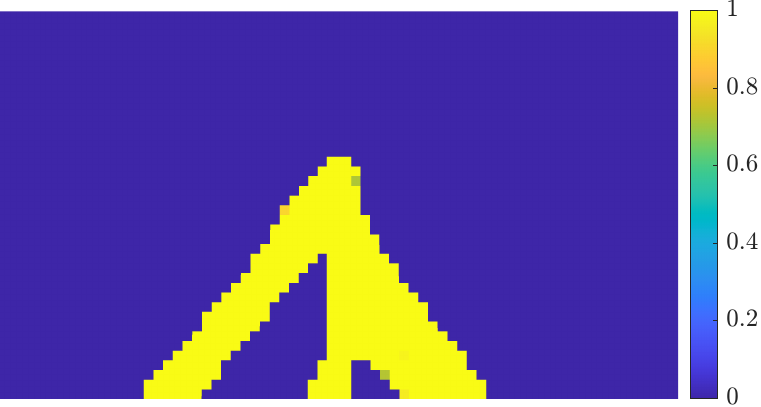}} \
	\subfloat[$p_\alpha = 18.0e-7$, $f = 0.152950$, \\ \vspace{1pt} $DM = 0.031043 \ \%$ \label{ch7_ffig:5k}]{\includegraphics[width=0.32\textwidth]{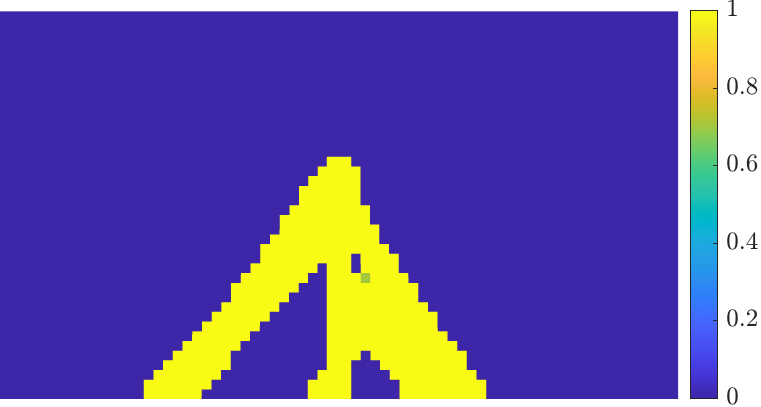}}
	\caption{Effect of $p_\alpha$ on the optimized designs. $\delta_{E|\mathit{\Upsilon}}$ is set to 1. $f$ is the objective function while $DM$ is the discreteness measure.}
	\label{ch7_fig:effect_of_p_alpha}
\end{figure}

\begin{figure}[t!]
	\centering
	\captionsetup[subfigure]{justification=centering}
	\subfloat[$p_\alpha = 8.0e-7$.]{\includegraphics[width=0.32\textwidth]{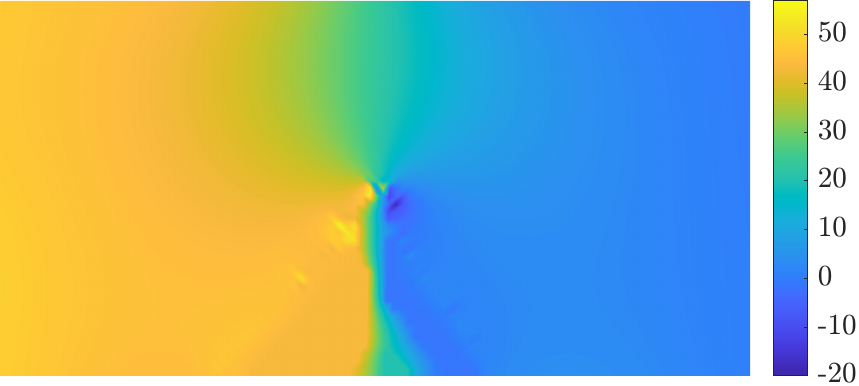}} \
	\subfloat[$p_\alpha = 8.5e-7$.]{\includegraphics[width=0.32\textwidth]{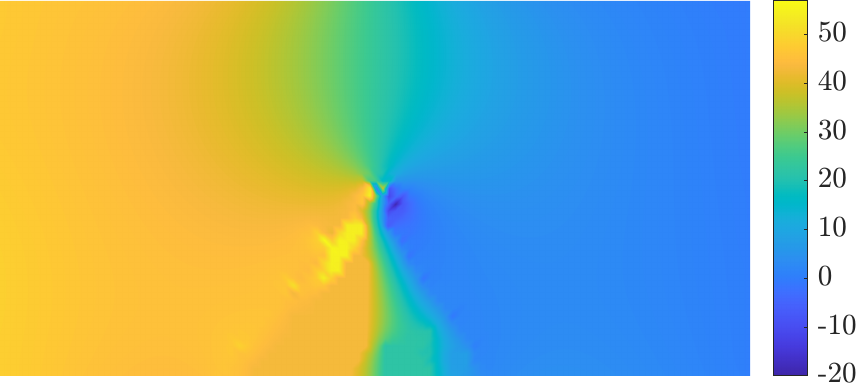}} \
	\subfloat[$p_\alpha = 9.0e-7$.]{\includegraphics[width=0.32\textwidth]{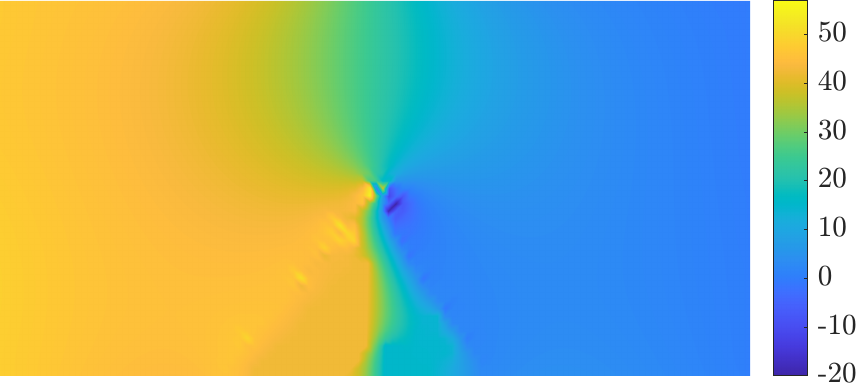}} \\
	\subfloat[$p_\alpha = 9.5e-7$.]{\includegraphics[width=0.32\textwidth]{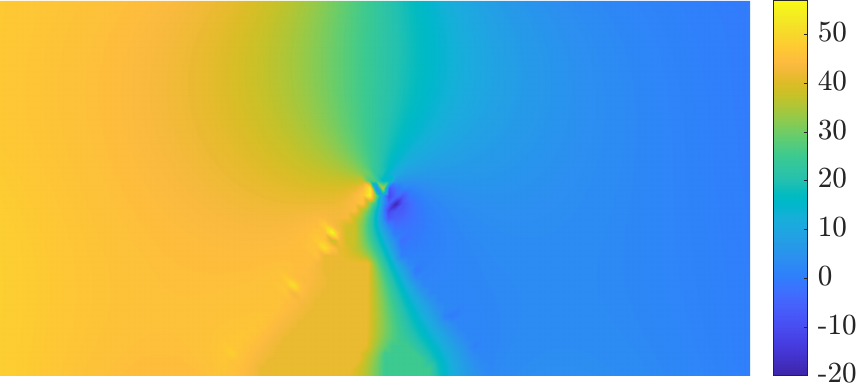}} \
	\subfloat[$p_\alpha = 10.0e-7$.]{\includegraphics[width=0.32\textwidth]{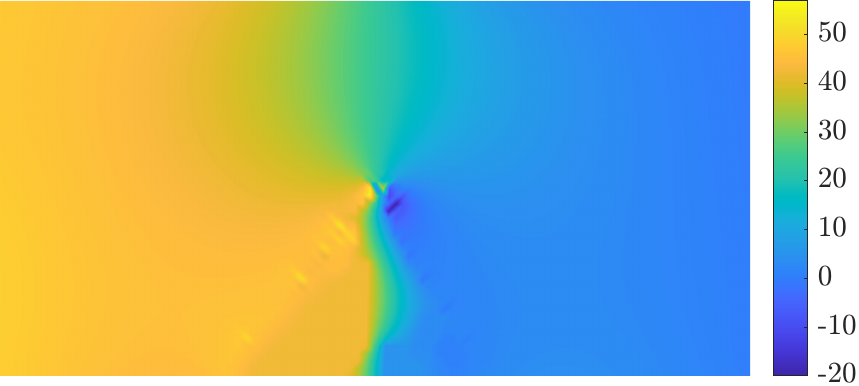}} \
	\subfloat[$p_\alpha = 11.5e-7$.]{\includegraphics[width=0.32\textwidth]{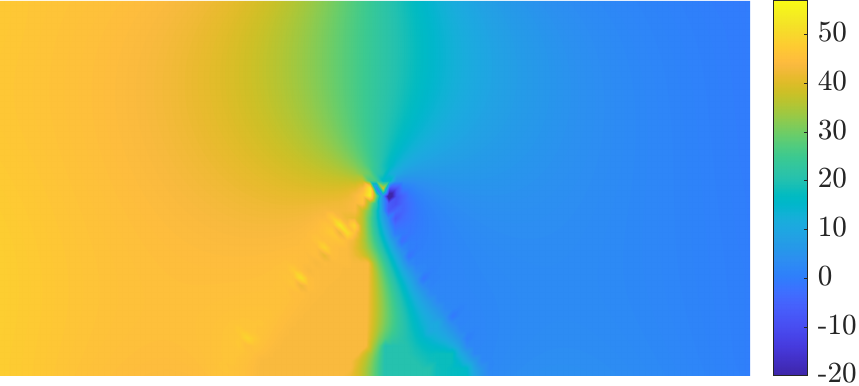}} \\
	\subfloat[$p_\alpha = 12.5e-7$.]{\includegraphics[width=0.32\textwidth]{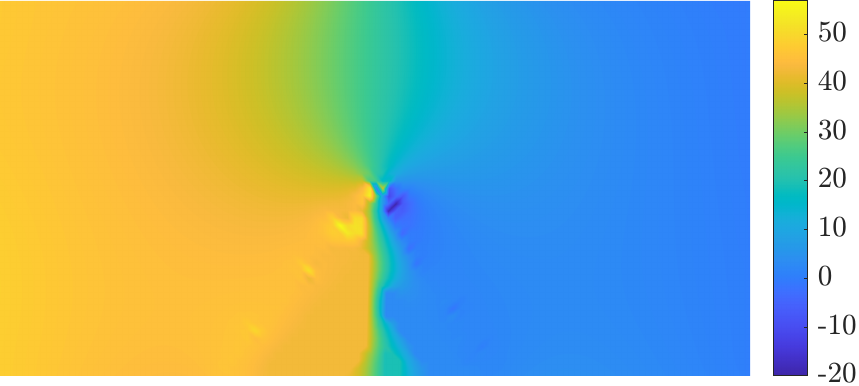}} \
	\subfloat[$p_\alpha = 14.0e-7$.]{\includegraphics[width=0.32\textwidth]{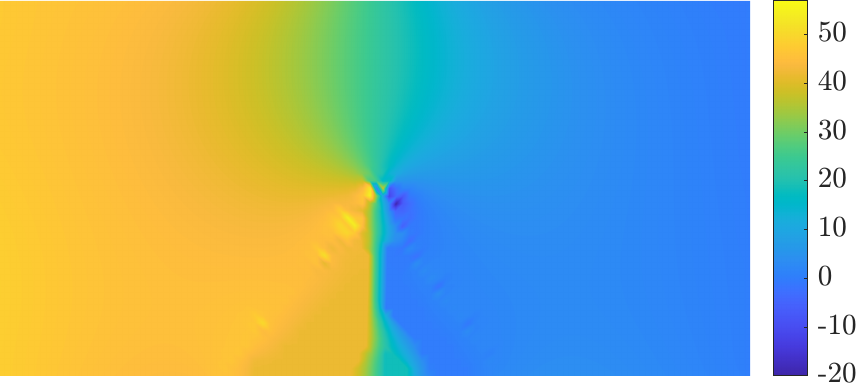}} \
	\subfloat[$p_\alpha = 15.5e-7$.]{\includegraphics[width=0.32\textwidth]{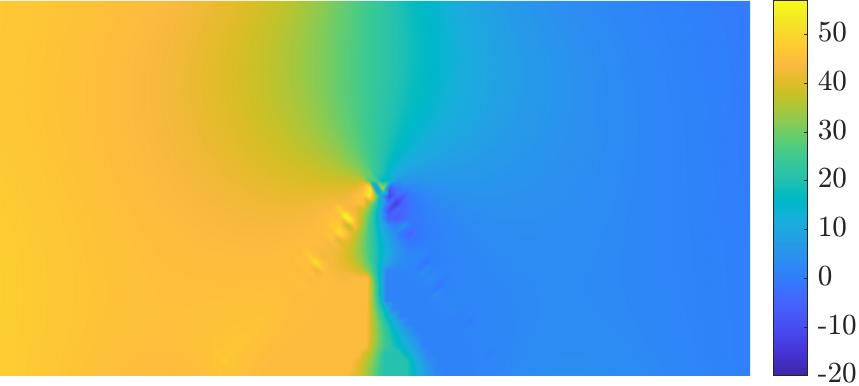}} \\
	\subfloat[$p_\alpha = 17.5e-7$.]{\includegraphics[width=0.32\textwidth]{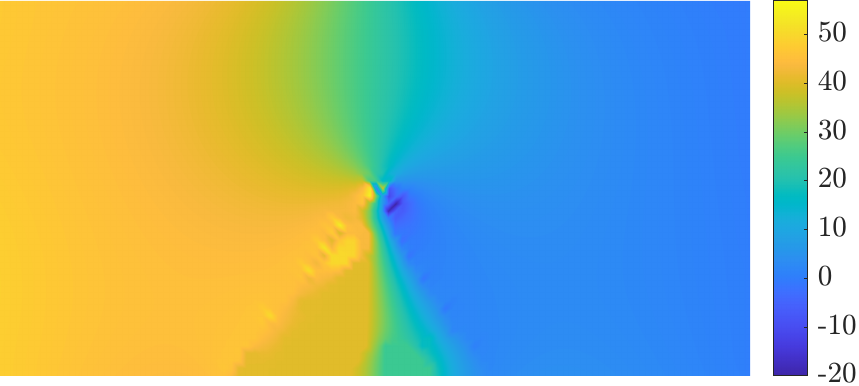}} \
	\subfloat[$p_\alpha = 18.0e-7$.]{\includegraphics[width=0.32\textwidth]{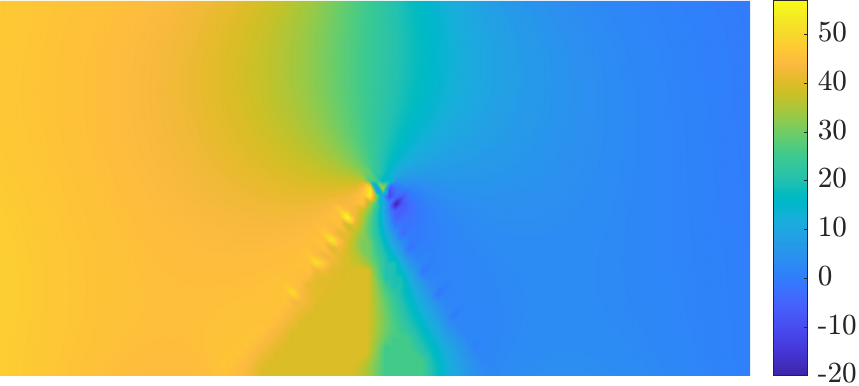}}
	\caption{Pressure fields for the optimized designs in Fig. \ref{ch7_fig:effect_of_p_alpha}.}
	\label{ch7_fig:prssr_flds_for_dfrnt_p_alpha}
\end{figure}

In a previous work \cite{Abdelhamid2022}, we investigated the effect of the projection and interpolation parameters on the optimized designs for TOFSI problems with fixed mesh. We summarize the main points here before discussing the appropriate parameters for including mesh deformation. In the absence of mesh deformation, the following was noted \citep[p.~13]{Abdelhamid2022}:
\begin{enumerate}[label=\alph*.]
	\item Implementing continuation on $p_\alpha$ doesn't have any convexification effects on the optimized designs, hence a single $p_\alpha$ value is selected for the entirety of the optimization process.
	\item The value of $p_\alpha$ is critical in that it determines the sensitivity of the fluid flow to intermediate density features, where pressure is far more sensitive to $p_\alpha$ than velocity.
	\item The proper value of $p_\alpha$ is the value that is large enough that the fluid flow can ``sense'' intermediate density features, yet small enough that the optimizer doesn't abuse these features to divert the flow away from the main structure in a non-physical manner.
	\item Penalizing the elastic modulus more than the force coupling tends to get rid of intermediate density features faster. In other words, a value of $p_E$ greater than $p_\alpha$ produces more discrete designs.
\end{enumerate}

For TOFSI problems with mesh deformation, we build upon these observations to determine the appropriate projection and interpolation parameters. One critical note is on how the existence of intermediate density features affects convergence towards discrete designs. In TOFSI problems with fixed mesh, an improper selection of $p_\alpha$ may result in the existence of intermediate density features in the final design without any impact on the convergence behavior of the problem (cf. \cite{Abdelhamid2022} for details). However, for TOFSI problems with mesh deformation, an improper $p_\alpha$ that produces intermediate density features tends to conflict with the raising of the robust formulation project parameter $\beta$. As $\beta$ is increased, these intermediate density features tend to become poorly supported or completely disconnected (in what is known as free-floating islands \cite{Lundgaard2018}), which ultimately results in failure of convergence. In our previous experiments on the column in a channel test problem in the absence of mesh deformation, the design tended to self-correct once these free-floating islands appeared, and no additional measures were needed. However, from our numerical experiments on TOFSI problems with mesh deformation, either the problem converges to a discrete design or it fails to converge at all. From our experience, the convergence failure occurs in the segregated analysis step. In the following, we discuss the selection of $p_\alpha$ then that of $\delta_{E|\mathit{\Upsilon}}$.

For the selection of $p_\alpha$, we experimented with a range of values while setting $\delta_{E|\mathit{\Upsilon}}$ to 1, i.e. $p_E$ and $p_\mathit{\Upsilon}$ experience the same degree of penalization. Figure \ref{ch7_fig:effect_of_p_alpha} shows optimized designs for a range of $p_\alpha$ values where the objective function $f$ and the discreteness measure\footnote{The discreteness measure $DM$ is calculated following Eq. 41 in \citep[p.~415]{Sigmund2007}.} $DM$ are reported for each design. We attempted to solve cases for $p_\alpha$ ranging from $5e-7$ to $100e-7$, however, most of the designs diverged somewhere in the segregated analysis step. The issue usually appears directly after a continuation update around $\beta = 16$ or $\beta = 32$ which is typically a strong projection with a steep descrease in the discreteness measure, thus favorable conditions for the appearance of poorly-supported features (or worse; free-floating islands). It can be clearly noted that in Fig. \ref{ch7_fig:effect_of_p_alpha}, all optimized designs that achieved convergence are almost free of intermediate density features as is evident in the very good discreteness measures reported under each design. There is little variation in both the objective functions and the discreteness measures reported in Fig. \ref{ch7_fig:effect_of_p_alpha}, however, a critical factor that helps in selecting the proper $p_\alpha$ is taking a close look at the pressure fields presented in Fig. \ref{ch7_fig:prssr_flds_for_dfrnt_p_alpha}. While the pressure drop across the solid, non-design column in the middle is very distinct, the pressure drop across the structural member upstream is not clear for the lower $p_\alpha$ values. This fact is also observed in our previous work on TOFSI with fixed mesh in \cite{Abdelhamid2022}, where the system tends to abuse the improper $p_\alpha$ by directing material to better support the non-design column yet lower their density values that the fluid doesn't ``see'' it, hence a lower objective function overall.

\begin{figure}[b!]
	\centering
	\captionsetup[subfigure]{justification=centering}
	\subfloat[$\delta_{E|\mathit{\Upsilon}} = 1/1.25$, $f = 0.154699$, \\ \vspace{1pt} $DM = 0.112550 \ \%$]{\includegraphics[width=0.32\textwidth]{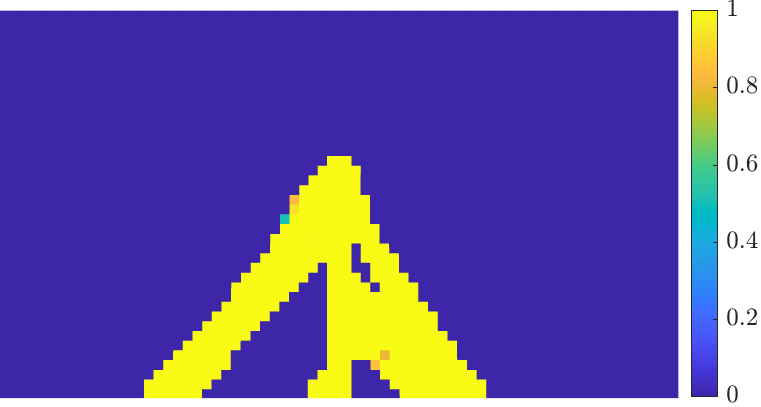}} \
	\subfloat[$\delta_{E|\mathit{\Upsilon}} = 1/1.50$, $f = 0.153183$, \\ \vspace{1pt} $DM = 0.050038 \ \%$]{\includegraphics[width=0.32\textwidth]{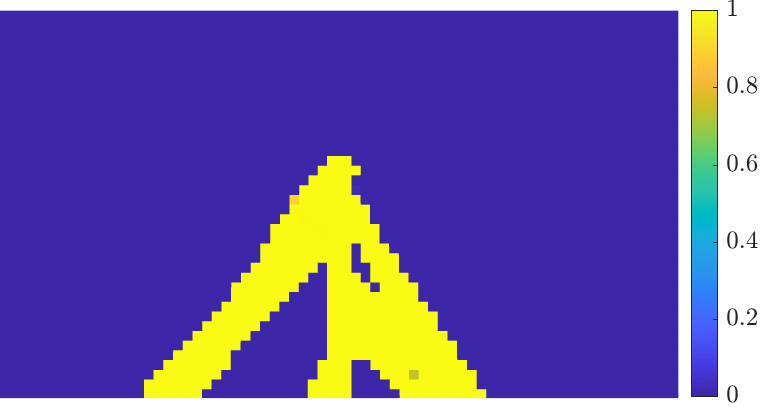}} \
	\caption{Effect of $\delta_{E|\mathit{\Upsilon}}$ on the optimized designs. $p_\alpha$ is set to $18e-7$.}
	\label{ch7_fig:pe_vs_pupsilon}
\end{figure}

For $\delta_{E|\mathit{\Upsilon}}$, we first attempted to use $\delta_{E|\mathit{\Upsilon}} > 1$ similar to previous work on TOFSI with fixed mesh \cite{Abdelhamid2022}, however, convergence couldn't be achieved for cases with $\delta_{E|\mathit{\Upsilon}} = 1.25$ and $\delta_{E|\mathit{\Upsilon}} = 1.50$. Hence, we had to utilize values lower than 1 for $\delta_{E|\mathit{\Upsilon}}$, such that $p_\mathit{\Upsilon}$ is more penalized than $p_E$. Figure \ref{ch7_fig:pe_vs_pupsilon} shows two optimized designs for $\delta_{E|\mathit{\Upsilon}} = 1/1.25$ and $\delta_{E|\mathit{\Upsilon}} = 1/1.50$. While the two designs don't show improvement in the objective function or the discreteness measure, the total solution time improves slightly as $\delta_{E|\mathit{\Upsilon}}$ is decreased. Nonetheless, further numerical experiments are needed to confirm this trend.

\subsection{Fixed vs Deforming Mesh}
\begin{figure}[b!]
	\centering
	\includegraphics[width=0.32\textwidth]{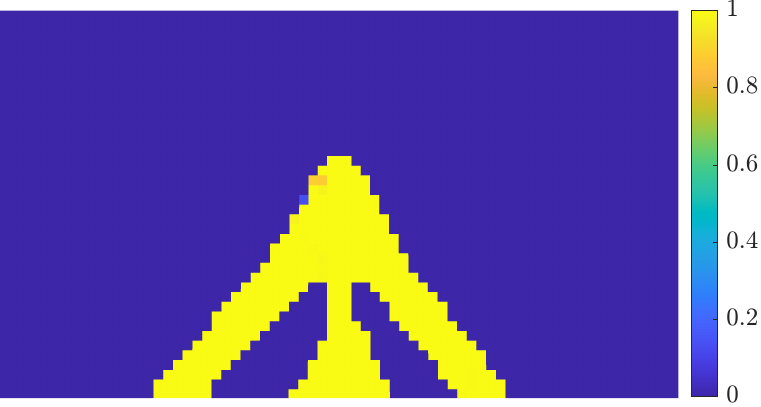}
	\caption{Optimized design for $p_\alpha = 18.0e-7$ without mesh deformation. $f = 0.159440$ and $DM = 0.060621 \ \%$.}
	\label{ch7_fig:fixed_mesh}
\end{figure}

The last numerical is to illustrate the benefits of implementing mesh deformation in TOFSI problems. For $p_\alpha = 18e-7$ and $\delta_{E|\mathit{\Upsilon}} = 1$, a design is generated while neglecting structural and mesh deformation in solving the fluid flow (see Fig. \ref{ch7_fig:fixed_mesh}). The optimized design in Fig. \ref{ch7_fig:fixed_mesh} is more symmetric and significantly different from those optimized under mesh deformation conditions. A cross-check is performed to differentiate the performance of each design under fixed and moving mesh conditions with the results reported in Table \ref{ch7_tab:fixed_vs_moving_mesh_crsschck}. The design optimized for moving mesh performs better when analyzed under both moving and fixed mesh conditions, which is a clear advantage for considering mesh deformation in TOFSI problems. However, the difference between the objective functions of both designs is small, hence a case maybe made for neglecting mesh deformations when the final deformation is typically small.

\begin{table}[t!]
	\centering
	\caption{Cross-check for the optimized designs with moving vs fixed mesh.}
	\label{ch7_tab:fixed_vs_moving_mesh_crsschck}
	\begin{tabular}{wc{0.7in} wc{0.7in} wc{1.25in} wc{1.25in}} \toprule
		& & \multicolumn{2}{c}{Obj. fun. analyzed for} \\
		Fig. & Designed for & \multicolumn{1}{c}{Moving Mesh} & \multicolumn{1}{c}{Fixed Mesh} \\ \midrule
		\ref{ch7_ffig:5k} & Moving Mesh & $0.152950$ & $0.157078$ \\ 
		\ref{ch7_fig:fixed_mesh} & Fixed Mesh & $0.156193$ & $0.159440$ \\ 
		\bottomrule
	\end{tabular}
\end{table}

\section{Conclusions}
\label{ch7_sec:conclusion}

In this work, we extend the application of density-based topology optimization to high-fidelity, fully-coupled fluid-structure interaction problems. We utilize the three-field formulation to describe the fluid-structure interaction problem where the arbitrary Lagrangian-Eulerian approach is used to deform the fluid mesh as a pseudo-structural system. The sensitivity analysis is derived using the discrete adjoint approach and verified numerically using the complex-step derivative approximation method. The main observations in this work are:
\begin{enumerate}
	\item The introduction a segregated analysis loop such that the fluid mesh reflects the structural deformations adds to the nonlinearity of the problem.
	
	\item Careful selection of the projection and interpolation parameters is critical to achieve convergence.
	
	\item Consideration of the structural deformations in the fluid mesh clearly has an advantage compared to one-way coupled fluid-structure interactions, albeit this advantage comes at a high computational cost. Hence, the compromise between the high-fidelity and computational costs is dependent on the problem of interest.
\end{enumerate}

\section*{Acknowledgments}

The authors are grateful to Prof. Krister Svanberg, KTH Royal Institute of Technology in Stockholm, for providing the MATLAB implementation of the Method of Moving Asymptotes. This research was enabled in part by support provided by BC DRI Group and the Digital Research Alliance of Canada (alliancecan.ca).

\bibliographystyle{elsarticle-num-names}
\pdfbookmark[0]{References}{References}
\bibliography{../../../library}

\end{document}